\documentclass[11pt]{article}
\usepackage{amssymb}
\usepackage{latexsym}

\usepackage{mathptmx}
\newtheorem{thm}{Theorem}[section]
\newtheorem{prop}[thm]{Proposition}
\newtheorem{lemma}[thm]{Lemma}

\newcommand{\gm }{\Gamma }

\newcommand{\reals}{\mathbb R}

\newcommand{\mathbold}{\bf}
\newcommand{\cinf}{C^{\infty}}

\def\qed{\rule{2.3mm}{2.3mm}}

\def\lbr {\lbrack\! \lbrack}
\def\rbr {\rbrack\! \rbrack}
\def\lfl {\lfloor\! \lfloor}
\def\rfl {\rfloor\! \rfloor}

\newcommand{\half}{\frac{1}{2}}

\begin{document}
\title{\bf Nambu-Dirac Structures on Lie Algebroids}

\author{
 A\"{\i}ssa WADE \\
 {\small Department of Mathematics,  The Pennsylvania State University} \\
  {\small University Park, PA 16802.}\\
   {\small {\it e-mail:} wade@math.psu.edu}
}

\date{}

\maketitle
\begin{abstract}
 The theory of Nambu-Poisson structures on manifolds
   is extended to the context of Lie algebroids
 in a natural way  based on the Vinogradov
 bracket associated with Lie algebroid cohomology.
 We show that, under certain assumptions, any 
 Nambu-Poisson structure on a Lie algebroid is decomposable. Also, we introduce the concept of a higher order Dirac structure on  a Lie algebroid.
 This allows to describe both Nambu-Poisson structures on Lie algebroids
  and Dirac structures on manifolds in the same setting. 

\end{abstract}

\noindent {\small{\bf Mathematics Subject Classification (2000):} 53Cxx, 58Axx,
 17B63, 81S10.}

\smallskip
\noindent {\small {\bf Key words:}  Leibniz algebras, Lie algebroid cohomology,
 Nambu-Poisson brackets, Dirac structures.}

\section {Introduction}

 The purpose of this Letter is to explore
  an unusual approach to Nambu-Poisson manifolds using
  the Vinogradov bracket defined below and the theory of Dirac structures.
  Our first attempt to describe Nambu-Poisson manifolds
 in terms of Dirac structures can be found in~\cite{W}.
 In the present paper, we work in a context which is general enough to include
 Nambu-Poisson manifolds and triangular Lie bialgebroids~\cite{MX}.

 Nambu-Poisson structures originate
  from an idea presented in~\cite{Na}, where Nambu took some tentative steps 
 toward a generalization of Hamiltonian mechanics. In 1975, Bayen and 
 Flato developed further Nambu's generalized mechanics (see~\cite{BF}).
  Takhtajan gave an axiomatic formalism for Nambu's bracket in~\cite{T}. 
   Since then, there has been a great deal of interest in the study
 of  Nambu-Poisson structures producing a beautiful series of papers
 by mathematicians and physicists.\par

 Consider  a Lie algebroid $E$ over a smooth manifold $M$. 
 It is well known that the Lie algebroid structure on $E$ induces a 
 differential operator $d_{\scriptscriptstyle E}$ on the space of all
 smooth sections of $ \Lambda^{\bullet}E^*$, where $E^*$
 is the dual of $E$. The Vinogradov bracket  is a derived bracket in 
 the sense of Kosmann-Schwarzbach (see~\cite{CV},~\cite{K-S}). 
   Precisely, the \textit{Vinogradov bracket} associated with
 $E$ is the $\reals$-bilinear operation on the space of
 graded endomorphisms of $\gm( \Lambda^{\bullet}E^*)$ given by 
 $\lbr a , b \rbr_{_E}= [[a, d_{\scriptscriptstyle E}], b]$, 
for any $a,\ b \in \mbox{End}(\Gamma(\Lambda^{\bullet} E^*))$,
where $[ \cdot,\cdot]$ is the usual commutator on endomorphisms of 
 $\gm(\Lambda^{\bullet} E^*)$. In Theorem~\ref{inj}, we show that there
 exists an injective map $\sigma$ from
$\gm (E \oplus \Lambda^{\bullet} E^*) $  to
 ${\rm End}(\Gamma(\Lambda^{\bullet} E^*))$  and
 an $\reals$-bilinear operation $\lfl \cdot  ,\cdot  \rfl_{_E}$
  on $\gm (E \oplus \Lambda^{\bullet} E^*) $ such that
\smallskip

\noindent (1) \   for any $e \in \gm (E \oplus \Lambda^{\bullet} E^*)$, 
 $\sigma (e)$ is $\cinf$-linear;\par

  \noindent (2) \  for any $e \in \gm (E \oplus \Lambda^{\bullet} E^*)$
 and for any $f \in \cinf(M)$, \ $\sigma(f e)=f \sigma(e)$;\par

 \noindent (3) \  for any $e_1, e_2 \in \gm (E \oplus \Lambda^{\bullet} E^*)$,
 \ $ \sigma(\lfl e_1  ,\ e_2 \rfl_{_E})=
 \lbr \sigma(e_1)  , \sigma (e_2) \rbr_{_E}$.\par

\noindent Precisely, $\sigma$ is defined by $\sigma(X, \alpha)=
 i_X + m(\alpha)$, where $i_X$ is the interior product by $X$ and
 $m(\alpha)$ is the exterior multiplication by $\alpha$. 
  The bracket $\lfl \cdot  ,\cdot  \rfl_{_E}$ is given by:
$$ \lfl (X,\alpha), \ (Y,\beta) \rfl_{_E}=
 ([X,Y]_{_ E}, \ {\cal L}^{^E}_X \beta -i_Y d_{\scriptscriptstyle E} \alpha),
 \leqno(1.1)$$

\noindent for any $(X,\alpha), \ (Y,\beta) \in
 \gm(E \oplus \Lambda^{\bullet} E^*)$,
 where ${\cal L}^{^E}_X =[ i_X, \ d_{\scriptscriptstyle E}]$.
  In view of this correspondence, one may think of 
$\gm (E \oplus \Lambda^{\bullet} E^*)$ as locally generated by 
 field operators satisfying the canonical commutation relations.
  The elements of $\gm(E)$ correspond to boson fields and the
 elements of $\gm(E^*)$ to fermion fields.  \par

 In the second part of this work,
 we exploit the bracket $\lfl \cdot  ,\cdot  \rfl_{_E}$, and in this way
 we are led to the definition of a
  Nambu-Poisson structure on a Lie algebroid $E$. A
 {\em Nambu-Poisson structure} of order $p$ on $E$ is  a smooth section
 $\Pi \in \gm(\Lambda^p E)$ such that  
$$[\Pi \alpha, \Pi]_{_E} \beta=-\Pi(i_{\Pi \beta}d_{\scriptscriptstyle E}
 \alpha), \ \mbox{for any} \ \alpha, \ \beta \in \gm(\Lambda^{p-1}E^*).$$ 

\noindent This definition includes both triangular Lie bialgebroids
 (see [MX]) and Nambu-Poisson structures on manifolds, as shown below.
  We prove that, given an arbitrary
 Nambu-Poisson structure $\Pi$ of order $p\geq 3$ on a Lie algebroid $E$
 with base $M$, $\Pi(d_{\scriptscriptstyle E} f)$
 is decomposable for
 any $f \in \cinf(M)$. In particular, $\Pi$ is decomposable when
 $\gm(E^*)$ is generated by the elements of the form 
 $d_{\scriptscriptstyle E} f$.

Finally, we study  Nambu-Poisson structures on a Lie algebroid $E$ from the 
 point of view of the theory of Dirac structures.  Very recently,
 the concept of a Nambu-Dirac manifold has been introduced and studied
 in~\cite{H}. This extension of  the definition of a Dirac manifold is 
 essentially based on  the existence of a Leibniz algebroid structure
 on the vector bundle  $\Lambda^{p-1} (T^*M)$ for a Poisson manifold $M$
 of order $p$. Note that there is another
  Leibniz algebroid structure on  $\Lambda^{p-1} (T^*M)$
  presented in~\cite{ILMP}.
  Here, inspired by the remarkable interplay between
  skew-symmetric bilinear forms and bivectors arising
  from linear Dirac structures~\cite{C},
   we use an alternative approach to
  define a Dirac structure of higher order
 on a Lie algebroid. We then generalize  both  strong Nambu-Dirac~\cite{H}
  and  Dirac structures on manifolds.
  We show that, to any 
 Nambu-Poisson structure $\Pi$ of order $p\geq 2$ on a Lie algebroid $E$,
 there corresponds  a Dirac structure of order $p$. Furthermore, we prove that
  any Dirac structure of order $p$ on $E$ gives rise to a Nambu-Poisson bracket
 on the set of its admissible functions. Thus,
  we recover certain results obtained in~\cite{H}. 
 Our main results are Theorem~\ref{inj},
 Theorem~\ref{decomp} , and Theorem~\ref{admis} \par

The paper is organized as follows. In Section 2,
 we review basic definitions and properties
 of Leibniz algebras and algebroids, as well as 
Lie algebroids. Then we define the Vinogradov bracket associated
 with a Lie algebroid and
  establish Theorem~\ref{inj}. 
 In Section 3, we define Nambu-Poisson structures on Lie algebroids
 and study their properties.
Finally, in Section 4, we study the relationships
 between Nambu-Poisson structures on Lie algebroids
  and the theory of Dirac structures.
 
\section{ Leibniz Algebras and Vinogradov Bracket}

\subsection{Definitions}
Throughout the text, all vector bundles have finite dimension.
We will be using Leibniz algebras and Lie algebroids in many parts of this
 work, so we briefly recall the definitions.\par

 A {\em left Leibniz bracket} 
 on a  real vector space  ${\cal V}$ (see~\cite{L}) 
  is given by  an $\reals$-bilinear operation $[\cdot  ,\cdot  ]_{_{\cal V}}: \
 {\cal V} \times
 {\cal V} \rightarrow {\cal V}$  such that
$$[a, [b,c]_{_{\cal V}}]_{_{\cal V}}=[[a,b]_{_{\cal V}}, c]_{_{\cal V}}+
[b, [a,c]_{_{\cal V}}]_{_{\cal V}}
,$$
\noindent for any   $a$, $b$, $c \in {\cal V}$.
  Then the pair $({\cal V}, [\cdot  ,\cdot  ]_{_{\cal V}})$ is called 
  {\em left Leibniz algebra}. For short, it will be called 
 Leibniz algebra.   The bracket
  $[\cdot  ,\cdot  ]_{_{\cal V}}$ becomes  a  Lie bracket if,
 additionally, it is skew-symmetric. \par

\smallskip
 Let $({\cal V}_1, [\cdot  ,\cdot  ]_1)$ and 
$ ({\cal V}_2, [\cdot  ,\cdot  ]_2)$ be two  Leibniz algebras.
 A morphism $\phi$ from $ \ ({\cal V}_1, [\cdot  ,\cdot  ]_1)$
  to $({\cal V}_2, [\cdot  ,\cdot  ]_2)$  is said
 to be a morphism of Leibniz algebras if
 $\phi([a ,b ]_1)= [\phi(a)  , \phi(b) ]_2,$
 for any $a$, $b \in {\cal V}_1.$ \par 

\smallskip
A {\em Leibniz algebroid} (see~\cite{ILMP}) consists of a vector bundle
$E \rightarrow M$ with a
  Leibniz bracket $[ \ , \ ]_{_E}$ on the space
  $\Gamma(E)$ of smooth sections of $E$ and a  bundle map 
 $\varrho :  E  \rightarrow TM$, extended to a map 
 between sections of these bundles such that the following properties are
 satisfied:
 $$[X_1, fX_2]_{_E}= f[X_1,X_2]_{_E}+ (\varrho(X_1)f)X_2
 \quad \mbox{and} \quad 
 \varrho([X_1, X_2]_{_E})=[\varrho(X_1),\varrho(X_2)],$$
\noindent for any smooth sections $X_1$, $X_2 \in \gm(E)$ 
 and for any smooth
  function $f$  defined on $M$. 
 Then $\varrho$ is called the {\em anchor} of the Leibniz algebroid. \par

 A {\em Lie algebroid} is a Leibniz algebroid $(E, [ \ , \ ]_{_E},
 \varrho)$  over $M$ such that $[ \ , \ ]_{_E}$ is skew-symmetric.

 \subsection{The Vinogradov bracket of a Lie algebroid}
 
Let $E$ be a Lie algebroid over a smooth $n$-dimensional manifold
 $M$ with anchor map $\varrho$ and let $E^*$ be its dual. Let
 $\gm(\Lambda^{\bullet} E^*)=\oplus_{k \geq 0} \Gamma(\Lambda^k E^*)$ 
  denote the exterior algebra generated by all
 smooth sections of $E^*$. Define
 $ d_{\scriptscriptstyle E}: ~~ \Gamma(\Lambda^{k-1}E^*) \longrightarrow
 \Gamma(\Lambda^k E^*)$ by:
 \begin{eqnarray*}
 \label{d_E}
 (d_{\scriptscriptstyle E} \xi)  (e_1, ..., e_k) & = &  \sum_i
 (-1)^{i+1} \varrho(e_i) (\xi (e_1, ..., \hat{e}_i, ..., e_k)) \\ &   + &
 \sum_{i < j} (-1)^{i+j} \xi([e_i, e_j]_{_E}, ..., \hat{e}_i, ...,
 \hat{e}_j, ..., e_k), \nonumber
 \end{eqnarray*}
 \noindent for any $\xi \in \Gamma(\Lambda^{k-1}E^*)$, and for any 
 $e_1, \dots , e_k \in \gm(E)$. It 
  satisfies $d_{\scriptscriptstyle E}^2 =
 0$. The cohomology of $(\Gamma(\Lambda^{\bullet} E^*),
 d_{\scriptscriptstyle E})$ is called the Lie algebroid cohomology of
 $E$ and denoted by $H^{\bullet}(E)$. 
Let $\mbox{End}^{\ell}(\Gamma(\Lambda^{\bullet} E^*))$  denote the space of
 all graded endomorphisms of degree $\ell$ on $\Gamma(\Lambda^{\bullet} E^*)$, and
$\mbox{End}( \Gamma(\Lambda^{\bullet} E^*))
= \oplus_{\scriptscriptstyle \ell \geq 0} \mbox{End}^{\ell}
 ( \Gamma(\Lambda^{\bullet} E^*))$.

\noindent Then, the {\em Vinogradov bracket} associated with the Lie algebroid
 $E$ is the $\reals$-bilinear operation
 $\lbr \cdot , \cdot \rbr_{_E}$ defined by 
$\lbr a , b \rbr_{_E}= [[a, d_{\scriptscriptstyle E}], b]$,
for any $a$,
 $b \in \mbox{End}(\Gamma(\Lambda^{\bullet} E^*))$,
   where $[ \cdot,
 \cdot]$ is the usual commutator on endomorphisms of 
 $\gm(\Lambda^{\bullet} E^*)$,
  i.e. $$[a, b]= a\circ b- (-1)^{|a| \ |b|} b \circ a,$$
\noindent when $a$, $b$ are graded endomorphisms  with
  respective degrees $|a|$, $|b|$.  \par
 For instance, the tangent bundle $TM$  of a smooth manifold $M$
 defines a Lie algebroid with respect to the usual Lie bracket
  on vector fields.
 The Vinogradov bracket associated with $TM$ was used in~\cite{CV} as a tool
 for extending an arbitrary Poisson bracket defined on $\cinf(M)$
  to the space of  co-exact forms. 
   As mentioned above, the Vinogradov bracket associated with
 a Lie algebroid $E$ described here
  is a particular case of the general construction of derived brackets 
  presented in~\cite{K-S}.
 In fact, what Kosmann-Schwarzbach called Vinogradov bracket
 in~\cite{K-S} is the skew-symmetrization of $\lbr \cdot , \cdot \rbr_{_E}$.
 But, we  will use the non skew-symmetric version in this work. 
  The main property of $\lbr \cdot , \cdot \rbr_{_E}$ (see~\cite{K-S}) is that
 $$ \lbr a, \lbr b, c\rbr_{_E} \rbr_{_E}=
   \lbr \lbr a,  b \rbr_{_E}, c\rbr_{_E} + (-1)^{(|a|+1)(|b|+1)}
   \lbr b, \lbr a, c\rbr_{_E} \rbr_{_E}, $$

\noindent for any graded endomorphisms  $a$, $b$, and $c$ with
  respective degrees $|a|$, $|b|$, and $|c|$ .

\smallskip  
Let $i_X$ denote the interior product by $X \in \gm(E)$ and 
 $m(\alpha)$ the multiplication by $\alpha \in \Gamma(\Lambda^{\bullet} E^*)$,
  i.e.
$m(\alpha) ( \eta)= \alpha \wedge \eta,
 \quad \forall \ \eta \in \Gamma(\Lambda^{\bullet} E^*)$.
 We associate the $\cinf(M)$-linear map 
 $\sigma: \  \ \gm(E \oplus \Lambda^{\bullet} E^*)
 \rightarrow   \mbox{End}(\Gamma(\Lambda^{\bullet} E^*))$ given by
$\sigma(X, \alpha)= i_X + m(\alpha)$,
for any $(X,\alpha) \in \gm(E \oplus \Lambda^{\bullet} E^*)$.
 With these notations, we have:

\begin{thm}
\label{inj}
 Let $(E, [ \cdot , \cdot]_{_E}, \varrho)$ be a Lie algebroid over $M$
 and  let  $\lbr \cdot , \cdot \rbr_{_E}$ be its corresponding Vinogradov
 bracket. Then, there is a  $\reals$-bilinear operation 
 $\lfl \cdot  ,\cdot  \rfl_{_E}$
 on $\gm(E \oplus \Lambda^{\bullet} E^*)$ such that 
 $\sigma: \  \ \gm(E \oplus \Lambda^{\bullet} E^*)) , 
 \  \lfl \cdot  ,\cdot  \rfl_{_E}) 
 \longrightarrow
   (\mbox{End}(\Gamma(\Lambda^{\bullet} E^*) , \ 
  \lbr \cdot  ,\cdot   \rbr_{_E})$
 is an injective morphism  of Leibniz algebras. Precisely,
  $\lfl \cdot  ,\cdot  \rfl_{_E}$ is defined by:
$$
 \lfl (X,\alpha), \ (Y,\beta) \rfl_{_E}=
 ([X,Y]_{_ E}, \ {\cal L}^{^E}_X \beta -i_Y d_{\scriptscriptstyle E} \alpha),
\leqno (1.1)$$

\noindent for any $(X,\alpha), \ (Y,\beta) \in
 \gm(E \oplus \Lambda^{\bullet} E^*)$,
 where ${\cal L}^{^E}_X =[ i_X, \ d_{\scriptscriptstyle E}]$.
\end{thm}

To prove  this theorem, we need the following lemma.
\begin{lemma}
\label{prep1}
Let $X \in \gm(E)$, $\alpha, \beta \in \Gamma(\Lambda^{\bullet} A^*)$. Then,
\begin{itemize}
 \item[]{ \rm (i)} \ \  $[i_X, \ m(\alpha)]= m(i_X \alpha)$;
 \item[]{ \rm (ii)} \ \ $[{\cal L}^{^E}_X, \  m(\alpha)]=
  m({\cal L}^{^E}_X \alpha)$;
\item[]{\rm (iii)} \ \ $[d_{\scriptscriptstyle E}, \ m(\alpha)]= m( d_{\scriptscriptstyle E}
 \alpha)$;
\item[]{\rm (iv)}  \ \ $[m(\alpha), \ m(\beta)]=0$;
\item[]{\rm (v)}  \ \ $ [{\cal L}^{^E}_X, {\cal L}^{^E}_Y]=
  {\cal L}^{^E}_{_{[X,Y]_{_E}}}.$
\end{itemize}
\end{lemma}

\noindent The proof of this lemma is straightforward.

\noindent{\it Proof of Theorem~\ref{inj}:} For any
 $X, Y \in \gm(E)$, we have
$\lbr i_X, i_Y \rbr_{_E}= [{\cal L}^{^E}_X, \ i_Y]=
i_{[X , Y]_{_E}}$.
Applying Property (ii) of Lemma~\ref{prep1}, we get 
$\lbr i_X, \ m(\beta) \rbr_{_E}= [{\cal L}^{^E}_X, \  m(\beta)]
= m({\cal L}^{^E}_X \beta)$,
for  any $X \in \gm(E)$,
 $\beta \in \Gamma(\Lambda^{\bullet} E^*)$.
  Using (iii), (i) and the symmetry property of the commutator,  one gets:
$$\lbr m(\alpha), \ i_Y \rbr_{_E}=
   [[m(\alpha), \ d_{\scriptscriptstyle E}], \ i_Y]
= -[ i_Y, m(d_{\scriptscriptstyle E} \alpha)]= 
 - m(i_Y d_{\scriptscriptstyle E} \alpha),$$
\noindent for  any $Y \in \gm(E)$, $\alpha \in \Gamma(\Lambda^{\bullet} A^*)$.
Finally, Property (iv) implies 
$\lbr m(\alpha), \ m(\beta) \rbr_{_E}=0,$
for  any $\alpha, \ \beta \in \Gamma(\Lambda^{\bullet} A^*)$.
 Thus combining the above terms, we obtain 
$$\lbr i_X + m(\alpha), \ i_Y + m(\beta) \rbr_{_E}=
 i_{[X , Y]_{_E}}+ m({\cal L}^{^E}_X \beta -  i_Y d_{\scriptscriptstyle E} \alpha),$$ 
\noindent for any $X, Y \in \gm(E)$, $\alpha, \beta \in 
\Gamma(\Lambda^{\bullet} E^*)$.
 Equivalently,

$$\lbr \sigma(X, \alpha), \  \sigma(Y,\beta) \rbr_{_E}=
\sigma ([X , Y]_{_E}, \ {\cal L}^{^E}_X \beta -  i_Y d_{\scriptscriptstyle E} \alpha)
= \sigma(\lfl (X, \alpha), \ (Y,\beta) \rfl_{_E}).$$

\noindent  Since $\sigma$ is injective and $\lbr \ , \ \rbr_{_E}$ is 
 a Leibniz bracket, it follows that $\lfl \cdot , \cdot \rfl_{_E}$ is
 a Leibniz bracket. Hence, $\sigma$ is an injective
  morphism of Leibniz algebras.
This completes the proof of Theorem~\ref{inj}.

\hfill \qed

\smallskip
\noindent{\bf Example 1.}  Let $E=TM$ be the
 tangent bundle of a smooth manifold $M$. Then, applying Theorem~\ref{inj}, we
 get  an operation $\lfl \cdot  ,\cdot  \rfl$ on $\gm(TM \oplus
  \Lambda^{\bullet}T^*M)$  which is ``isomorphic'' to the Vinogradov bracket
 via the map $\sigma$. Namely,
  $\lfl \cdot  ,\cdot  \rfl$  is given by

$$\lfl (X, \alpha), \ (Y, \beta) \rfl = ([X, Y], \ L_{_X} \beta -i_Y
  d \alpha),$$

\noindent where $[X, Y]$ is the Lie bracket of the vector fields $X$ and $Y$,
  $L_X$ is the Lie derivative with respect to 
 the vector field $X$ and $d$ is the 
de Rham differential. The restriction of $\lfl \cdot  ,\cdot  \rfl$  to
$\gm(TM \oplus T^*M)$ coincides with the Courant bracket (see~\cite{C}).

\smallskip
\noindent{\bf Remark.} 
 Theorem~\ref{inj} pushes further an
  observation made by Kosmann-Schwarzbach in an unpublished paper, where
  it is shown that the Courant bracket is obtained by using 
  $\sigma_{|\gm(TM \oplus T^*M)}$.

\section{ Nambu-Poisson Structures}

\subsection{Nambu-Poisson Structures on Manifolds}
Let $M$ be a  smooth $n$-dimensional manifold.
A {\em Nambu-Poisson structure} on $M$ of order $p$ (with $2\leq p \leq n$) 
 is given by a $p$-vector field  $\Pi$ which satisfies the fundamental 
 identity:
\begin{eqnarray*}
 \{f_1,...,f_{p-1}, \{g_1,...,g_p\}_{_\Pi} \}_{_\Pi} =
 \sum_{k=1}^{p}  \{g_1,...,g_{k-1}, \{f_1,...,f_{p-1},g_k\}_{_\Pi},
 g_{k+1},..., g_p\}_{_\Pi} 
\end{eqnarray*}

\noindent  for any $ f_1,\dots ,f_{p-1}, g_1, \dots g_p \in C^{\infty}(M)$,
  where $\{ \ \ \}_{_\Pi}$ is  defined by 

$$\{f_1, \dots ,f_p \}_{_\Pi}=\Pi(df_1, \dots , df_p) \quad \forall
  f_1,\dots ,f_{p} \in C^{\infty}(M).$$ 

 \noindent A manifold equipped with such a structure is called 
 {\em Nambu-Poisson manifold}.  Nambu-Poisson  structures of order 2 On $M$
 are  Poisson structures.  It is easy to see
 that the fundamental identity  is equivalent  to
 $[X_{f_1, \dots , f_{p-1}}, \Pi]_{_{SN}}=0, \quad \mbox{for any } \ 
f_1,\dots ,f_{p-1} \in C^{\infty}(M),$
 where $[\cdot , \cdot]_{_{SN}}$
 is the Schouten-Nijenhuis on multi-vector
 fields and  $X_{f_1, \dots , f_{p-1}}$ is the Hamiltonian
 vector field given by
$\langle X_{f_1, \dots , f_{p-1}}, \ dg \rangle
=\{f_1, \dots , f_{p-1}, g\}_{_\Pi}.$

\subsection{Nambu-Poisson structures on Lie algebroids}

It is  known that, given a  a Lie algebroid
   $(E, [ \cdot , \cdot]_{_E}, \varrho)$ over $M$, the algebra
 $\gm(\Lambda^{\bullet}E)$ endowed with the exterior product and the
 generalized Schouten bracket determines a Gerstenhaber algebra.
  Consider a smooth section $\Pi \in \gm(\Lambda^p E)$. We say that $\Pi$ is a
 {\em Nambu-Poisson structure} of order $p$ on $E$ if 
$$[\Pi \alpha, \Pi]_{_E} \beta=-\Pi(i_{\Pi \beta}d_{\scriptscriptstyle E}
 \alpha), \ \mbox{for any} \ \alpha, \ \beta \in \gm(\Lambda^{p-1}E^*).$$ 

In particular, this property implies that 
 $[\Pi \alpha, \Pi]=0$ when $ d_{\scriptscriptstyle E}\alpha=0$.
 This definition agrees with that of a Nambu-Poisson structure on
 a manifold, which is obtained when $E=TM$. \par
 
It is  natural to ask if a Nambu-Poisson  structure of order $p>2$ on 
 a Lie algebroid $E$ is decomposable. It is well know that such
 a phenomenon occurs
 for  Nambu-Poisson structures of order $p>2$ on a
 smooth manifold $M$. It turns out that  
 any  Nambu-Poisson structure of order $p>2$ on 
 a Lie algebroid $E$  is decomposable if 
   $\gm(E^*)$  is generated by  the sections $d_{\scriptscriptstyle E}f$,
 for $f \in \cinf(M)$.
  Precisely, we have the following result:

\begin{thm}
 \label{decomp}
 Let $\Pi$ be a 
 Nambu-Poisson  structure of order $p \geq 3$ on 
 a Lie algebroid $E$. For any $f \in \cinf(M)$, 
 $\Pi(d_{\scriptscriptstyle E}f)$ is decomposable. In particular,
 if $\gm(E^*)$ is generated by elements of the form 
 $d_{\scriptscriptstyle E}f$ then $\Pi$ is decomposable. 
\end{thm}

\noindent{\it Proof}: Assume that $\Pi$ is 
 a Nambu-Poisson  structure of order $p \geq 3$ on $E$. 
 Take $f \in \cinf(M)$, 
 $\eta \in \gm(\Lambda^{p-1}E^*)$ and set $\alpha =
 f d_{\scriptscriptstyle E}f \wedge \eta$.
On the one hand, for any $\beta \in \gm(\Lambda^{p-2}E^*)$ we have
\begin{eqnarray*}
[\Pi \alpha, \ \Pi]_{_E} \beta =-\Pi(i_{\Pi \beta}d_{\scriptscriptstyle E}
 \alpha)= f \Pi\Big(i_{\Pi \beta}(d_{\scriptscriptstyle E}f 
 \wedge d_{\scriptscriptstyle E} \eta)\Big) .
\end{eqnarray*}

\noindent On the other hand,
\begin{eqnarray*}
[\Pi \alpha, \ \Pi]_{_E} \beta &=&
 [f \Pi(d_{\scriptscriptstyle E}f \wedge \eta), \ \Pi]_{_E} \beta \cr
 &=& -\Big( \pm \Pi(d_{\scriptscriptstyle E}f) \wedge
 \Pi (d_{\scriptscriptstyle E}f \wedge \eta) +
 f[ \Pi(d_{\scriptscriptstyle E}f \wedge \eta), \ \Pi]_{_E} \Big) \beta \cr
 &=& -\Big( \pm\Pi(d_{\scriptscriptstyle E}f) \wedge
  \Pi (d_{\scriptscriptstyle E}f
  \wedge \eta)\Big)\beta +
 f \Pi\Big(i_{\Pi \beta}(d_{\scriptscriptstyle E}f 
 \wedge d_{\scriptscriptstyle E} \eta)\Big).
\end{eqnarray*}

\noindent Therefore, 
$$\Pi(d_{\scriptscriptstyle E}f) \wedge \Pi (d_{\scriptscriptstyle E}f
 \wedge \eta) =0, \
 \mbox{for any} \ f \in \cinf(M), \ \eta \in \gm(\Lambda^{p-2}E^*).$$

\noindent Hence, $\Pi(d_{\scriptscriptstyle E}f)$ is decomposable
  (see~\cite{MVV}).
 There follows the result.

\hfill \qed

\begin{prop}
\label{Poisson}
 Consider a Lie algebroid
 $E$ with base $M$.  Then $\Pi$ is a 
 Nambu-Poisson structure of order 2 on $E$
 if and only if $[\Pi, \Pi]_{_E}=0$.
\end{prop}

 \noindent This proposition shows that the definition of
 Nambu-Poisson structure on $E$ agrees with that of a Poisson structure on $E$
 when $p=2$.
  To  prove it, we need the following lemma
 which has been shown in~\cite{K-SM} and~\cite{LX}.

\begin{lemma}
\label{formula}
  For any $\eta_1$ and $\eta_2 \in \gm(E^*)$, we have
 $$\Pi\Big({\cal L}^{^E}_{_{\Pi\eta_1}} \eta_2 -
  {\cal L}^{^E}_{_{\Pi\eta_2}} \eta_1
 -  d_{\scriptscriptstyle E}( \Pi(\eta_1, \eta_2)) \Big) -
  [ \Pi\eta_1, \Pi\eta_2]_{_E} =
 \half [\Pi, \Pi]_{_E} (\eta_1 \wedge \eta_2).$$
\end{lemma}
 
\noindent{\it Proof of Proposition~\ref{Poisson}:} 
 We have $[ \Pi\eta_1, \Pi\eta_2]_{_E}=
 [ \Pi\eta_1, \Pi]\eta_2 + \Pi({\cal L}^{^E}_{_{\Pi\eta_1}} \eta_2)$.
Add  $\Pi(i_{_{\Pi\eta_1}} d_{\scriptscriptstyle E}\eta_2)$ in both sides
 of this last equation. We get
$$
[ \Pi\eta_1, \Pi\eta_2]_{_E}-
\Pi({\cal L}^{^E}_{_{\Pi\eta_1}} \eta_2) +
 \Pi(i_{_{\Pi\eta_1}} d_{\scriptscriptstyle E}\eta_2) =
   [ \Pi\eta_1, \Pi]_{_E} \eta_2 +
  \Pi(i_{_{\Pi\eta_2}} d_{\scriptscriptstyle E}\eta_1) .$$
 \noindent Now using Lemma~\ref{formula}, we get
$\half [\Pi, \Pi]_{_E} (\eta_1 \wedge  \eta_2) =-
  \Big([\Pi\eta_1, \Pi]_{_E} \eta_2+
  \Pi(i_{_{\Pi\eta_2}} d_{\scriptscriptstyle E}\eta_1)\Big)$.
Hence,
$[\Pi, \Pi]_{_E} =0 \iff [\Pi\eta_1, \Pi]_{_E} \eta_2 +
  \Pi(i_{_{\Pi\eta_2}} d_{\scriptscriptstyle E}\eta_1)=0, $
 for any $\eta_1$ and $\eta_2 \in \gm(E^*)$. % This completes the proof.

\hfill \qed

\smallskip

Now we are going to show that any 
 Nambu-Poisson structure of order $p$ on a Lie algebroid
 $E$ with base $M$ gives rise to a natural Nambu-Poisson bracket of order $p$
 on the algebra $\cinf(M)$ of functions on $M$. Define
 the bracket $\{....\}_{_\Pi}$ by
$$\{f_1, \dots , f_p \}_{_\Pi}=\Pi(d_{\scriptscriptstyle E} (f_1)\wedge 
 \dots \wedge  d_{\scriptscriptstyle E}(f_p )),$$
\noindent for $f_1, \dots , f_p  \in \cinf(M)$. Obviously, this bracket 
 is skew-symmetric and
 $\{f_1, \dots , f_{p-1}, \cdot \}_{_\Pi}$ is a derivation
 with respect to the pointwise product of functions. We denote by 
$X_{f_1, \dots , f_{p-1}}$ the corresponding vector field. Then, we get
$${\cal L}^{^E}_{X_{f_1, \dots , f_{p-1}}} g=\{f_1, \dots , f_{p-1}, g\}_{_\Pi}.$$
\noindent Furthermore, since
  $[\Pi(d_{\scriptscriptstyle E} (f_1)\wedge  \dots \wedge 
  d_{\scriptscriptstyle E}(f_{p-1})), \ \Pi]=0$, it follows
\begin{eqnarray*}
\{f_1, \dots , f_{p-1}, \{g_1, \dots , g_p \}_{_\Pi} \}_{_\Pi}&=&
 {\cal L}^{^E}_{X_{f_1, \dots , f_{p-1}}} \{g_1, \dots , g_p \}_{_\Pi} \cr
 &=& \sum_{i=1}^p \{g_1, \dots , g_{i-1}, \{f_1, \dots , f_{p-1}, g_i\}_{_\Pi},  \dots, g_p \}_{_\Pi}.
\end{eqnarray*}

\noindent Therefore, we get the following result:

\begin{prop}  Given a Lie algebroid
   $(E, [ \cdot , \cdot]_{_E}, \varrho)$ with base $M$,
 there is a natural Nambu-Poisson bracket on $\cinf(M)$ arising
 from any Nambu-Poisson structure on $E$.
\end{prop}

\begin{prop}
\label{algebroid}
 Given   a Lie algebroid
   $(E, [ \cdot , \cdot]_{_E}, \varrho)$ over $M$
   and a Nambu-Poisson structure $\Pi$ on $E$, 
the triplet $(\Lambda^{p-1}E^*, \lfl \cdot, \cdot \rfl_*, \varrho \circ \Pi)$
determines a Leibniz algebroid, where $\lfl \cdot, \cdot \rfl_*$ is defined
 by 
$$\lfl \alpha, \beta \rfl_* = {\cal L}^{^E}_ {\Pi \alpha}
  \beta -  i_{\Pi \beta} d_{\scriptscriptstyle E} \alpha,
\ \mbox{for any} \ \alpha, \ \beta \in \gm(\Lambda^{p-1}E^*) $$
\end{prop}

\noindent{\it Proof:} By simple computations,
 one gets
 $$\lfl \alpha, \ f \beta \rfl_*= f\lfl \alpha, \ \beta \rfl_*
+ (\varrho(\Pi \alpha)f) \beta \quad \mbox{and}
 \quad \varrho (\Pi\lfl \alpha, \beta \rfl_*)  =
  [ \varrho (\Pi\alpha), \ \varrho (\Pi\beta) ],$$

\noindent for any $\alpha, \ \beta \in \gm(\Lambda^{p-1}E^*)$. 
 Moreover, we have: 
\begin{eqnarray*}
\lfl \alpha, \lfl \beta , \gamma \rfl_* \rfl_*&=&
{\cal L}^{^E}_{\Pi\alpha}(\lfl \beta, \gamma \rfl_*) 
-i_{\Pi\lfl \beta, \gamma \rfl_*} d\alpha \cr
&=&{\cal L}^{^E}_{\Pi\alpha}(L_{\Pi\beta} \gamma) -
{\cal L}^{^E}_{\Pi\alpha}(i_{\Pi\gamma} d\beta)\cr &&
- {\cal L}^{^E}_{\Pi\beta}(i_{\Pi\gamma} d\alpha)
 + i_{\Pi\gamma}(L_{\Pi\beta} d\alpha).
\end{eqnarray*}

\noindent By similar computations we obtain 
\begin{eqnarray*}
\lfl \lfl \alpha, \beta \rfl_* , \gamma \rfl_*+
\lfl \beta , \lfl \alpha, \gamma \rfl_* \rfl_*&=&
{\cal L}^{^E}_{\Pi\alpha}({\cal L}^{^E}_{\Pi\beta} \gamma) -
{\cal L}^{^E}_{\Pi\alpha}(i_{\Pi\gamma} d\beta) \cr &&
- {\cal L}^{^E}_{\Pi\beta}(i_{\Pi\gamma} d\alpha)
 + i_{\Pi\gamma}({\cal L}^{^E}_{\Pi\beta} d\alpha).
\end{eqnarray*}

\noindent Thus,
$\lfl \alpha, \lfl \beta , \gamma \rfl_* \rfl_*=
\lfl \lfl \alpha, \beta \rfl_* , \gamma \rfl_*+
\lfl \beta , \lfl \alpha, \gamma \rfl_* \rfl_*$.
The proposition follows immediately.\par

\hfill \qed

\section{Higher order Dirac structures}

\subsection{Linear Dirac structures}
Let $V$ be an $n$-dimensional vector space over $\reals$ and $V^*$ its dual.
 Consider the canonical symmetric bilinear form $\langle \cdot , \cdot \rangle$ on $V \oplus V^*$ given by

$$\langle (v_1,  \omega_1), \ (v_2,  \omega_2) \rangle=
\half (\omega_1(v_2)+ \omega_2(v_1)).$$

A {\em linear Dirac structure} on $V$ (see~\cite{C}) is a subspace $L$ of $V \oplus V^*$
which is maximally isotropic with respect to $\langle \cdot , \cdot \rangle$. 
There are two equivalent characterizations of linear Dirac structures
 given by the following propositions.
 
\begin{prop}{\rm ~\cite{C}}
There is a bijective correspondence between linear Dirac structures on 
$V$ and the pairs $(P, \Omega)$ consisting of a subspace $P$ of $V$ and a 
  skew-symmetric bilinear form $\Omega: P \times P \rightarrow \reals$
\end{prop}
\noindent{\it Proof:} Let $\pi : V \oplus V^* \rightarrow V$ and 
$\pi^* : V \oplus V^* \rightarrow V^*$ be the canonical projections. 
If $L \subset V \oplus V^*$ is a linear Dirac structure on $V$, then
 there exists a well defined skew-symmetric bilinear map
$\Omega_{L}: \pi(L) \times \pi(L) \rightarrow \pi(L)$ given by
 $\Omega_{L}(v_1, v_2)=\omega_1(v_2)$,
for any $(v_1,  \omega_1)$, $(v_2,  \omega_2) \in L$.

Conversely, let $P$ be a subspace of $V$ equipped with a 
skew-symmetric bilinear form $\Omega : P \rightarrow P^*$. Define
$L=\{(v, \omega) \in V \oplus V^* \ | \ v \in P \ {\rm and} \ 
\Omega(v)= \omega_{|P}\}$. $L$ is isotropic with respect to
   $\langle \cdot , \cdot \rangle$ since $\Omega$ is skew-symmetric.
  Moreover, 
$$\mbox{dim} L= \mbox{dim}P+ \mbox{dim ker} \pi_{|P}=
\mbox{dim} P+ \mbox{dim Ann}(P)=n,$$
  where ${\rm Ann}(P)$ denotes the annihilator
 of the subspace $P$. Therefore, $L$ is a linear Dirac structure on $V$.
\hfill \qed

\noindent Similarly, one gets the proposition:
\begin{prop}{\rm ~\cite{C}}
There is a bijective correspondence between linear Dirac structures on 
$V$ and  pairs $(S, \Pi)$ consisting of a subspace $S$ of $V^*$ and a 
  2-vector $\Pi \in \Lambda^2S^* $.
\end{prop}

\subsection{Higher order Dirac structures  on vector spaces}
All vector spaces considered in this section have finite dimension.

\smallskip
\noindent {\bf Definition.} 
Let $V$ and $W$ be two vector spaces.
 A {\em linear  quasi Dirac structure of order $p$
 on $V$ relative to $W$} consists of a subspace $P$ of $V$ together
 with a  skew-symmetric $p$-form $\Omega: \
  P \times \dots \times P \rightarrow W$. \

\smallskip
 For instance, let $V=W=\reals^n$. Consider a skew-symmetric 
 bilinear operation  $B$ on $\reals^n$. Then $(\reals^n, B)$ is a 
 linear  quasi Dirac structure of order 2
 on $\reals^n$ relative to $\reals^n$.
\medskip

\noindent{\bf Definition.}
A {\em linear Dirac structure of order $p$} on a real
 vector space $V$   relative to $\reals$
 is a  linear quasi Dirac structure $(P, \Omega)$ 
  of order $p$ on $V$ such that for any 
 $Z_1, \dots , Z_{p-1} \in \Lambda^{p-1}P$, there exists  $v \in P$ satisfying
  $$(i_{Z_1} \Omega) \wedge \dots \wedge (i_{Z_{p-1}}\Omega) =i_{v} \Omega.
 \leqno({\cal H})$$

\noindent Obviously, (${\cal H}$) is always satisfied  when $p=2$.
In other words, any linear quasi Dirac structure $(P, \Omega)$ 
  of order 2 is a linear Dirac structure of order 2. \par

 For a better understanding of the condition (${\cal H}$) for $p>2$, let
 us  consider the subspace of $P^*$ given by
$E_{\Omega}= {\rm Span} \ \{ i_{v_1} \dots i_{v_{p-1}} \Omega \ | \ 
 v_i \in P \}$. It is known that ${\rm dim} (E_{\Omega})  \geq p$ and 
 ${\rm dim} (E_{\Omega})=p$ if and only if $\Omega$ is decomposable.
 Let $\{ \eta_1, \dots , \eta_{r}\}$ be a basis for $E_{\Omega}$. Then,
$$\Omega= \sum_{1 \leq i_1 < \dots < i_{p-1} \leq r}
\eta_{i_1} \wedge \dots \wedge \eta_{i_{p-1}} \wedge \nu_{i_1\dots i_{p-1}},$$

\noindent where $\nu_{i_1\dots i_{p-1}} \in 
 E_{\Omega}$. The condition (${\cal H}$) means that none of the 1-forms
 $\nu_{i_1\dots i_{p-1}}$  is zero. We deduce that (${\cal H}$) is always
 true when $\Omega$ is decomposable.\par

 Now, consider  a linear Dirac structure $(P, \Omega)$ 
  of order $p$ on a real
 vector space $V$ relative
to $\reals$. On the one hand, $\Omega$ can be viewed as a map from
 $P$ to $\Lambda^{p-1} P^*$. In this case, we associate
 the subspace $L$  of $V \times \Lambda^{p-1} V^*$:
$$L=\{(v, \omega) \ | \
  v \in P \ {\rm and} \ 
i_{v}\Omega= \omega_{| \Lambda^{p-1}P} \}.$$

\noindent On the other hand,  regarding $\Omega$ as a map from
$\Lambda^{p-1} P$ to $P^*$, the natural subspace
 $L^{\perp} \subset \Lambda^{p-1} V \times  V^*$ to be defined is:

$$L^{\perp}=\{(Z, \eta)  \ | \
  Z \in \Lambda^{p-1} P \ {\rm and} \ 
i_{Z}\Omega= \eta_{| P} \}.$$

\noindent For any $(v, \omega) \in L$ and for any $(Z, \eta) \in L^{\perp}$,
  we have
$$\eta_{_{| P}}(v)+ (-1)^p \omega_{_{| \Lambda^{p-1}P}} (Z)=0,
\leqno (1.2)
$$
\noindent since $\Omega(Z)(v)= (-1)^{p-1}\Omega(v)(Z)$. 

\smallskip

\noindent{\bf Notation.} Let $\pi^*(L^{\perp}) $ denote
 the projection of $L^{\perp}$ onto 
 $V^*$ and  $\pi^*(L) $  the projection of $L$ onto 
 $\Lambda^{p-1}V^*$.  Then we have the following lemmas.

\begin{lemma}
  $${\rm Ann}(\pi^*(L^{\perp}))=L \cap V.$$
\end{lemma}

\noindent{\it Proof:} In fact, $v \in  {\rm Ann}(\pi^*(L^{\perp}))$
if and only if 
 $\eta (v)=0, \ \forall (Z, \eta) \in L^{\perp}.$
This is equivalent to say that
$\Omega(Z)(v)= (-1)^{p-1}\Omega(v)(Z) =0, \ \forall  Z \in \Lambda^{p-1}P$.
But,  $\Omega(v)=0$ means that $v \in L\cap V$. There follows the lemma. 

\hfill \qed

\begin{lemma}
$$ \Lambda^{p-1}\Big(\pi^*(L^{\perp}) \Big) \subset  \pi^*(L).$$
\end{lemma}

\noindent{\it Proof:} For any $(Z_1, \eta_1) \dots (Z_{p-1}, \eta_{p-1})
\in L^{\perp}$, we have
$$\eta_1 \wedge \dots \wedge \eta_{p-1} \ {_{|\Lambda^{p-1}P}}=
 (i_{Z_1} \Omega) \wedge \dots \wedge (i_{Z_{p-1}}\Omega).$$

\noindent We deduce from the definition 
 of a linear Dirac structure  of order $p$ the existence of 
 an element $v \in P$ such that
$i_v\Omega=\eta_1 \wedge \dots \wedge \eta_{p-1}\ {_{|\Lambda^{p-1}P}}$.
Thus  $(v, \eta_1 \wedge \dots \wedge \eta_{p-1})$
 is in $L$. There follows the lemma.

\hfill \qed

\smallskip 
Now, we set $S= \pi^*(L^{\perp})$ and identify $S^*$ with 
 $V/ L\cap V$. Consider the $p$-vector 
 $ \Pi  \in \Lambda^{p} S^*$ defined  by
$\Pi(\eta)= Z_{| \Lambda^{p-1}S}$, for any $(Z, \eta) \in L^{\perp}$.
 $\Pi$ is well defined. Indeed, if
 $(Z_1, \ \eta) $ and $(Z_2, \ \eta)$ are in $L^{\perp}$,
  then $(Z_1-Z_2, 0) \in L^{\perp}$. Let $\eta_1, \dots , \eta_{p-1}$
 be in $S$. There exists $v$ in $P$ such that 
$(v, \eta_1 \wedge \dots  \wedge\eta_{p-1}) \in L$. 
Thus, $(Z_1-Z_2)(\eta_1 \wedge \dots  \wedge\eta_{p-1})=0$.
This shows that $Z_1=Z_2$ on $\Lambda^{p-1}S$. 
 Since $\Pi$ is skew-symmetric, we get
$$ Z_1(\eta_2)+Z_2(\eta_1)_{| \Lambda^{p-2}S}=0,
  \ \forall \ (Z_1,\eta_1), \ (Z_2,\eta_2) \in L^{\perp}.$$

\noindent If $(v, \eta_1 \wedge \dots  \wedge\eta_{p-1}) \in L$,
 with  $\eta_i \in S$, then   
$\Pi(\eta_1 \wedge \dots  \wedge\eta_{p-1})=v_{|S}$.  
This follows from the fact that,
 for any $(Z_p, \eta_p) \in  L^{\perp}$, we have
\begin{eqnarray*}
(\Pi(\eta_1 \wedge \dots  \wedge\eta_{p-1}))(\eta_p) &=&
 (-1)^{p-1} (\Pi (\eta_p))(\eta_1 \wedge \dots  \wedge\eta_{p-1}) \cr
&=& (-1)^{p-1} Z_p (\eta_1 \wedge \dots  \wedge\eta_{p-1})=
 v(\eta_p).
\end{eqnarray*}

\noindent Furthermore, $\Pi$ has the following property: for any
 $\omega_1 , \dots , \omega_{p-1} \in \Lambda^{p-1} S$, there exists 
$\eta \in S$  such that 
$$\Pi(\omega_1) \wedge \dots \wedge \Pi(\omega_{p-1})= \Pi(\eta).
 \leqno({\cal H}')$$

In a similar way, one can show that, given any pair
 $(S, \Pi)$ of a subspace $S$ of $V^*$ and a $p$-vector
 on $S$ that satisfies $({\cal H}')$,
  there corresponds a linear Dirac structure $(P, \Omega)$ of order $p$ on 
 $V$. There follows the result:
\begin{prop}
\label{corresp}
There is a one-to-one correspondence between linear
 Dirac structures of order $p$ on a real vector space $V$
  relative to $\reals$ (with
 $2\leq p \leq dim V$) and the pairs $(S, \Pi)$ consisting
of a subspace $S \subset V^*$ and a $p$-vector
 $\Pi \in \Lambda^p S^* $ satisfying the property
{\rm (${\cal H}'$)}.
\end{prop}

\subsection{Higher order Dirac structures on Lie algebroids}

\noindent {\bf Definition.} 
 Let  $(E, [ \cdot , \cdot]_{_E}, \varrho)$ be a Lie algebroid over $M$.
 {\em A Dirac structure  of order $p$ on $E$} is a 
sub-bundle $L$ of $E\oplus \Lambda^p E^*$  which determines a linear
 Dirac structure of order $p$ relative to $\reals$
 at each point and whose sections are
  closed under the Leibniz bracket $\lfl \cdot , \cdot \rfl_{_E}$ defined 
 as in (1.1). The fact that $L$ determines a linear
 Dirac structure of order $p$ means that there exist a sub-bundle
 $P$ of $E$ and a skew-symmetric  bilinear form 
 $\Omega_{L}$ on $P$ satisfying (${\cal H}$) and such that
 the fibre of $L$ at a point $x \in M$ is
 $$L_x=\Big\{ (X, \alpha)_x  
 \ | \ X_x \in P_x \ {\rm and} \ (\Omega_{L}(X))_x=
  \alpha_{x |_{\Lambda^{(p-1)}P_x}}\Big \}.$$
 
\noindent Observe that if the sections of
 a sub-bundle $L$ of  $E\oplus \Lambda^p E^*$ are closed under 
   $\lfl \cdot , \cdot \rfl_{_E}$  then 
 $(L, \lfl \cdot , \cdot \rfl_{_E}, \varrho_{|L})$
  determines a Leibniz algebroid, where $\varrho_{|L}(X, \alpha)=
 \varrho(X)$ at every point.

\smallskip
\noindent{\bf Example 2.}
 A Dirac structure of order 2 on $TM$ is just a Dirac structure
 on $M$ in the sense of Courant (see~\cite{C}).\par

\smallskip
\noindent {\bf Example 3.} 
Let $\Pi$ be a 
 Nambu-Poisson  structure of order $p \geq 2$ on 
 a Lie algebroid $E$. With the notation 
 $X_{f_1 \dots f_{p-1}}=\Pi(d_{_E}f_1 \wedge\dots \wedge d_{_E}f_{p-1})$,
   consider
$$L={\rm Span} \ \{ (X_{f_1 \dots f_{p-1}}, \
  d_{_E}f_1 \wedge\dots \wedge d_{_E}f_{p-1}) \ | \ f_i \in \cinf(M)\}. $$

\noindent For $f_1, \dots , f_{p-1}, g_1, \dots , g_{p-1} \in \cinf(M)$, we 
 have
\begin{eqnarray*}
\lfl (X_{f_1 \dots f_{p-1}}, \ d_{_E}f_1 \wedge\dots \wedge d_{_E}f_{p-1}), \ 
 (X_{g_1 \dots g_{p-1}}, \ d_{_E}g_1 \wedge\dots \wedge d_{_E}g_{p-1})\rfl_{_E}
 = \cr
 \Big([X_{f_1 \dots f_{p-1}}, X_{g_1 \dots g_{p-1}}]_{_E}, \
\sum_{i=1}^{p-1} 
d_{_E}g_1 \wedge \dots   \wedge d_{_E}\{f_1,\dots, f_{p-1}, g_i\}_{_\Pi}
 \wedge \dots  \wedge d_{_E}g_{p-1} \Big).
\end{eqnarray*}
\noindent Moreover for any $(X, \alpha)$ and $(Y,\beta) \in \gm(L)$,
$f, g \in \cinf(M)$,

$$\lfl f(X, \alpha), \ g(Y,\beta)\rfl_{_E}= (fg [X, Y]_{_E}+
 f\langle d_{_E} g, X \rangle Y
 - g\langle d_{_E}f, Y \rangle X ,  \ \gamma),$$

\noindent where
$$\gamma= f g({\cal L}^{^E}_ {X}
  \beta -  i_{Y} d_{\scriptscriptstyle E} \alpha)+
 f\langle d_{_E} g , X \rangle \beta -
 g\langle d_{_E}f, Y \rangle \alpha+ g d_{_E}f\wedge(
 i_{X} \beta +i_{Y}\alpha).$$

\noindent But, one deduces from the fundamental identity
  that 
$$\Pi (d_{_E}f\wedge( i_{X} \beta +i_{Y}\alpha)) =0, $$

\noindent for any $X=X_{f_1 \dots f_{p-1}}$, \
 $\alpha =d_{_E}f_1 \wedge\dots \wedge d_{_E}f_{p-1},$ \
 $Y=X_{g_1 \dots g_{p-1}}$, and $\beta=d_{_E}g_1 \wedge\dots \wedge d_{_E}g_{p-1}$.  Therefore,
the sections of $L$ are closed the bracket $\lfl \cdot, \cdot \rfl_{_E}$. 
  So, using Proposition~\ref{corresp}, one sees that 
 $L$ determines a Nambu-Poisson structure of order $p$ on $E$. 

\smallskip

\noindent{\bf Example 4.}  Let $M$ be an $n$-dimensional
 manifold. 
   Recall from [H] that 
  an {\em almost Nambu-Dirac structure} of order $p$  on $M$ is a sub-bundle
 $L \subset  TM \times \Lambda^{p-1}T^*M$ such that
$\omega_1(X_2)+\omega_2(X_1)_{|\Lambda^{p-2}  \pi(L)}=0$,
for any $(X_1, \omega_1), \ (X_2, \omega_2) \in L$ and
$\Lambda^{p-1}\pi(L) =\pi( \widetilde{L})$, where 
$\widetilde{L}=\{(Z, \eta)  \in \Lambda^{p-1}TM \times T^*M \ | \
  \omega (Z)+(-1)^p \eta(X)=0, \ \forall  
 (X, \omega) \in L \}$, here $ \pi(L)$ and $\pi( \widetilde{L})$ denote
 the projections from $L$ and $\widetilde{L}$ onto their 
first  components, respectively. If the sections of $L$ are closed
 under the bracket defined as in (1.1), then $L$ is called a
 {\em Nambu-Dirac structure} of order $p$  on $M$.
 It is called a {\em strong Nambu-Dirac structure} of order $p$  on $M$
 if additionally 
$\Lambda^{p-1}(\pi^*( \widetilde{L})) \subset \pi^*(L)$,
 where $ \pi^*(L)$ and $\pi^*( \widetilde{L})$ denote
 the projections from $L$ and $\widetilde{L}$ onto their 
 second  components, respectively.
 It has been shown in
 [H] that an almost Nambu-Dirac structure of order $p$  on $M$
 induces a skew-symmetric $p$-form $\Omega_L$ on  $ \pi(L)$.
 Any  strong Nambu-Dirac structure
 of order $p$  on $M$ is  a
 Nambu-Dirac structure  of order $p$ on the Lie algebroid
  $E=TM$ as defined above.

\subsection{ The Nambu-Poisson Algebra of Admissible Functions}
 In this section, we give a definition of admissible functions for a 
Dirac structure of order $p$ on a Lie algebroid. Also, we 
show that, given any arbitrary
   Dirac structure of order $p$ on
 a Lie algebroid, there exists a Nambu-Poisson bracket defined on the set of
 admissible functions. \par

Let $(E, [ \cdot , \cdot]_{_E}, \varrho)$ be a Lie algebroid over $M$.
  Consider a sub-bundle $L$ of $E \oplus \Lambda^{p-2}E^*$.
 Assume that $L$ determines a Dirac structure of order $p$ on $E$, that is,
 $\gm(L)$ is closed under the bracket $\lfl \cdot  ,\cdot  \rfl_{_E}$
 defined as in (1.1) and there exist a sub-bundle
 $P$ of $E$ and a skew-symmetric 2-form 
 $\Omega_{L}$ on $P$ satisfying $({\cal H})$ such that
 the fibre of $L$ at a point $x \in M$ is
$L_x=\Big\{ (X, \alpha)_x  
 \ | \ X_x \in P_x \ {\rm and} \ (\Omega_{L}(X))_x=
  \alpha_{x |_{\Lambda^{p-1}P_x}}\Big \}$. 
A function $f \in \cinf(M)$ is said to be {\em $L$-admissible} if 
 $d_{\scriptscriptstyle E} f _{|L \cap E}=0.$

\noindent Now, we denote by $L^{\perp}$ the sub-bundle of
  $\Lambda^{p-1}E \oplus E^*$ whose
 fibre at a point $x \in M$ is given by:
 $$L^{\perp}_x= \Big\{ ({\mathbold Z}, \eta)_x
  \ | \ (\eta(X)+ (-1)^p \alpha({\mathbold Z}))_x=0, \  \forall \ (X, \alpha)_x
   \in L_x \Big\}.$$

In view of the above discussion, $f$ is $L$-admissible if and only if there
 exists a $(p-1)$-vector field $X_f$ such that 
 $e_f=({\mathbold X}_f, d_{\scriptscriptstyle E} f) \in \gm( L^{\perp})$. 
 Remark that for a given $L$-admissible function $f$,  $e_f$ is unique up to
 sections of $L^{\perp} \cap \Lambda^{p-2}E$. Let
  $f_1, \dots, f_p \in \cinf(M)$
 be $L$-admissible functions. Define
 $$\{f_1, \dots , f_{p} \}= {\mathbold X}_{f_1}(d_{\scriptscriptstyle E}f_2
  \wedge \dots \wedge d_{\scriptscriptstyle E}f_p )
= - {\mathbold X}_{f_2}(d_{\scriptscriptstyle E}f_1
  \wedge \dots \wedge d_{\scriptscriptstyle E}f_p )$$
\noindent This skew-symmetric bracket is well defined and is local.
 The product 
$gh$ of two $L$-admissible functions  $g$ and $h$ is $L$-admissible since
 $(g {\mathbold X}_h + h{\mathbold X}_g, d_{\scriptscriptstyle E}(gh))$
 becomes a section of $L^{\perp}$
 when $e_g=({\mathbold X}_g, d_{\scriptscriptstyle E} g)$ and
  $e_h=({\mathbold X}_h, d_{\scriptscriptstyle E} h)$ are
 sections of $\pi(L^{\perp})$. We have,
$$\{f_1, \dots , f_{p-1}, gh\}=
 g\{f_1, \dots , f_{p}, h\}+ h \{f_1, \dots , f_{p}, g\}.$$
 
\noindent Furthermore, we have the following result:

\begin{thm}
\label{admis}
 Let $L$ be a Dirac structure of order $p$ on a Lie algebroid $E$.
 Then the set of $L$-admissible functions has an induced 
 Nambu-Poisson bracket.
\end{thm} 

\noindent {\it Proof:} We only have to prove that the fundamental
 identity holds. But, this is an consequence of the fact that  
 the sections of $L$ are closed under the bracket 
$\lfl \cdot, \cdot \rfl_{_E}$.
 Indeed, if $f_1, \dots, f_{p-1}, g_1, \dots, g_{p-1}$ are  
  $L$-admissible  then 
  $(X_{f_1, \dots, f_{p-1}},  \ d_{\scriptscriptstyle E} f_1 \wedge
   \dots \wedge d_{\scriptscriptstyle E} f_{p-1})$ and 
$(X_{g_1, \dots, g_{p-1}}, d_{\scriptscriptstyle E} g_1 \wedge \dots \wedge
  d_{\scriptscriptstyle E} g_{p-1})$ are sections of $L$, where
$$X_{f_1, \dots, f_{p-1}}= X_{f_1}(d_{\scriptscriptstyle E} f_2\wedge
 \dots \wedge
  d_{\scriptscriptstyle E} f_{p-1}).$$
 \noindent Hence,  we have
$$\Big( \lfl X_{f_1, \dots, f_{p-1}}, \ X_{g_1, \dots, g_{p-1}}\rfl_{_E},
 \ {\cal L}_{X_{f_1, \dots, f_{p-1}}}^{^E}
(d_{\scriptscriptstyle E} g_1 \wedge \dots \wedge
  d_{\scriptscriptstyle E} g_{p-1}) \Big) \in \gm(L).$$

\noindent Since 
$${\cal L}_{X_{f_1, \dots, f_{p-1}}}^{^E}
\Big(d_{\scriptscriptstyle E} g_1 \wedge \dots \wedge
  d_{\scriptscriptstyle E} g_{p-1}\Big) = \sum_{i=1}^{p-1}
 d_{\scriptscriptstyle E} g_1 \wedge \dots \wedge 
d_{\scriptscriptstyle E} \{f_1, \dots, f_{p-1}, g_i\} \wedge \dots 
 \wedge d_{\scriptscriptstyle E} g_{p-1}.$$

\noindent We deduce from Equation (1.2) that for any arbitrary
 $L$-admissible function $h$, we have
$$ \Big\langle d_{\scriptscriptstyle E} h,  \lfl X_{f_1, \dots, f_{p-1}}, \ X_{g_1, \dots, g_{p-1}}\rfl_{_E}\Big \rangle=
 \sum_{i=1}^{p-1} \{g_1 , \dots ,
 \{f_1, \dots, f_{p-1}, g_i\} , \dots , g_{p-1}, h\}.$$
\noindent There follows immediately the the fundamental identity.
\hfill \qed

\smallskip
If a sub-bundle $L  \subset  E \times \Lambda^{p-1}E^*$
 determines a Dirac structure of order $p$ on $E$, then 
 $\varrho(L \cap E)$ is an integrable and singular distribution.
  Let ${\cal F}$ denote the singular foliation associated with
  $\varrho(L \cap E)$,
 which is called characteristic foliation of $L$. 
  Any $L$-admissible function is constant along the leaves of ${\cal F}$. 
 Assume that  the characteristic foliation is {\em simple}, that is, 
  ${\cal F}$ is a regular foliation such that $M/{\cal F}$ is a smooth manifold
 and ${\rm pr}: M \rightarrow M/{\cal F}$  is a submersion.
  In this case, $L$ is said to be {\it reducible}. 
The set of all $L$-admissible functions can be identified with 
 $\cinf(M/{\cal F})$. Thus, $M/{\cal F}$ has a Nambu-Poisson bracket. 
 It would be interesting
 to get a result analogous to Theorem 3.3 of~\cite{LWX} by showing that,
  given an integrable 
 distribution  whose foliation ${\cal F}$ is simple,
 any Nambu-Poisson structure of order $p$ on $ M/{\cal F}$ gives rise to
    a reducible Dirac structure of order $p$ on $E$. \par 
 We end this paper by noting that  the problem  of how to place  Nambu-Jacobi
 manifolds in the setting of Dirac structures is left open. 
  We hope that we will return to this question in the future.

\noindent {\bf Acknowledgement.}  I would  like to thank 
  J.-C. Marrero and D. Sternheimer for providing me with references
   and their helpful comments.
  Thanks go also to J. Stasheff for  many stimulating discussions
 about the Vinogradov bracket during Spring 2001.

\end{document}